\documentclass{amsart}
\usepackage{amssymb,amsmath,amsthm,latexsym,verbatim}

\issueinfo{00}
{01}
{08}
{2005}

\copyrightinfo{2005}{American Mathematical Society}

\newcounter{thmcount}
\newtheorem{theorem}[thmcount]{Theorem}
\newtheorem{corollary}[thmcount]{Corollary}
\newtheorem{definition}{Definition}[section]
\newtheorem{proposition}{Proposition}[section]
\newtheorem{lemma}{Lemma}[section]
\newtheorem*{remark}{Remark}

\newtheorem*{notation}{Notation}
\theoremstyle{definition}
\newtheorem{paragr}{}[section]
\begin{document}

\bibliographystyle{amsalpha}

\title{Mixing on Rank-One Transformations}

\author{Darren Creutz}
\address{Department of Mathematics, University of California: Los Angeles, 520 Portola Plaza, Los Angeles, CA 90095, USA}
\email{dcreutz@math.ucla.edu}

\author{Cesar E. Silva}
\address{Department of Mathematics, Williams College, Williamstown, MA 01267, USA} 
\email{csilva@williams.edu}


\date{February 3, 2006}

\begin{abstract}
We prove that mixing on rank-one transformations is equivalent to the 
spacer sequence being slice-ergodic.  Slice-ergodicity, introduced in 
this paper, generalizes the notion of ergodic sequence to the uniform 
convergence of ergodic averages (as in the mean ergodic theorem) over 
subsequences of partial sums.  We show that polynomial staircase 
transformations satisfy this condition and therefore are mixing.
\end{abstract}

\maketitle

\section{Introduction}
\begin{paragr}\textbf{Rank-One Transformations.} Rank-one 
transformations 
are transformations \lq\lq well-approximated" by a sequence of 
discrete spectrum  transformations, so it was very surprising when in 
1970  Ornstein \cite{Or72} showed the existence of rank-one {\it mixing} 
transformations.   Rank-one mixing  transformations are mixing of all 
orders \cite{Ka84}, \cite{Ry93} and enjoy other 
remarkable properties, see e.g. \cite{Ki88}.  Ornstein's construction 
was stochastic 
in nature:  there is a class of rank-one transformations 
so that almost surely a transformation in that class is mixing; however, 
it 
did not yield a deterministic procedure for constructing one.
\end{paragr}
\begin{paragr}\textbf{Staircase Transformations.} A few years later,   
Smorodinsky conjectured that a specific rank-one transformation, the
classic staircase transformation, is mixing.    In 1992, Adams and Friedman 
\cite{AF92} gave a deterministic algorithm involving a sequence of 
cutting and 
stacking constructions that produced  a mixing rank-one 
transformation, and later Adams \cite{Ad98} proved  that 
 Smorodinsky's conjecture is true.  Informally, a staircase transformation is a cutting and stacking transformation with  sequence $\{r_n\}$ of natural numbers such that at the $n^{\text{th}}$ stage the 
  $n^{\text{th}}$ column or stack is cut into $r_n$ subcolumns and 
\lq\lq spacers" (see Section \ref{S:rankone}) are placed in a staircase 
fashion 
on the subcolumns before stacking, i.e., the number of spacers in each 
subsequent subcolumn is increased  by $1$.  Adams showed that the 
resulting staircase transformation is mixing provided that $\frac 
{r_n^2}{h_n}\to 0$ 
as $n\to \infty$ (which also implies that $T$ is finite 
measure-preserving), where $h_n$ denotes the number of levels, or 
height,  of the  $n^{\text{th}}$ column.  
He then asked whether the mixing property 
holds for every finite measure-preserving staircase transformation simply 
under the assumption that $r_n\to\infty$.  
In 2003, Ryzhikov wrote the authors a short email stating
that in 2000 he  gave a lecture where
he  proved that all staircases are mixing  \cite{Ry03} (giving a
positive answer to Adams'  question), but no argument was included and
no preprint has been available. After this paper was completed, 
Ryzhikov informed the authors that his paper was forthcoming.
We would also like to thank Ryzhikov  for asking a question
that clarified our writing of Section~\ref{S:MixRankOne}.
The application of our main theorem shows that polynomial staircase 
transformations are mixing (Theorem~\ref{T:polyR1mix}).  Specializing to the case
of linear polynomials gives 
 another proof that staircase transformations are mixing.
\end{paragr}
\begin{paragr}\textbf{Restricted Growth.}
The $\frac{r_{n}^{2}}{h_{n}} \to 0$ condition, a restriction on the 
asymptotic growth of the spacers relative to the column height, was 
generalized to all rank-one transformations and 
called \lq\lq restricted growth" in \cite{CS04}.  The staircase 
transformation of Smorodinsky's conjecture is obtained when $r_n= n+1$; 
verifying that it satisfies the restricted 
growth condition is straightforward.  In \cite{CS04}, the authors   
proved an equivalence between mixing and a condition on the spacer 
sequence for rank-one transformations with restricted growth.  It followed that
 restricted growth rank-one transformations with the 
sequence of spacers given by a polynomial satisfying some general 
conditions (including the staircases of \cite{Ad98}) are mixing.  
Ornstein's result also follows from that theorem.
\end{paragr}
\begin{paragr}\textbf{Our Result.} In this 
paper we lift the restricted growth condition from the theorems in \cite{CS04}.   We 
introduce the notion of a slice-ergodic sequence and prove in 
Theorem~\ref{T:main} that a rank-one transformation is mixing if and 
only if 
its spacer sequence is slice-ergodic.  We use this theorem to show that 
all polynomial staircase transformations are mixing.
%Recently the authors have 
%obtained a version Theorem~\ref{T:main} for a larger
%class of transformations.
\end{paragr}

\section{Mixing Properties}\label{S:mixingprops}

\begin{paragr}\textbf{Dynamical Systems.}
For our study, \textbf{dynamical system} shall mean a 
standard probability 
measure space 
$(X,\mathcal{B},\mu)$ and \textbf{transformation} $T : X \to X$ that is 
invertible, measurable and measure-preserving.  Throughout the paper,
$X = [0,1)$, $\mu$ is Lebesgue measure on $X$ 
and $\mathcal{B}$ is the algebra of $\mu$-measurable subsets of $X$.
\end{paragr}

\begin{paragr}\textbf{Mixing.}
A transformation $T$ is \textbf{mixing}
when for all $A,B\in\mathcal{B}$,
\[
\lim_{n\to\infty}\mu(T^{n}(A) \cap B) - \mu(A)\mu(B) = 0;
\]
$\{t_{n}\}$ is a \textbf{mixing sequence} (with respect to 
$T$) when for all $A,B\in\mathcal{B}$,
\[
\lim_{n\to\infty}\mu(T^{t_{n}}(A) \cap B) - \mu(A)\mu(B) = 0.
\]
\end{paragr}

\begin{paragr}\textbf{Ergodicity.}
A transformation $T$ is
\textbf{ergodic} when for all $A\in\mathcal{B}$, if $T^{-1}(A) = A$ then 
$\mu(A) = 0$ or $\mu(A) = 1$.  The \textbf{mean (von Neumann) ergodic 
theorem} states that $T$ is ergodic if and only if for 
all $B\in\mathcal{B}$,
\[
\lim_{n\to\infty}\int\big{|}\frac{1}{n}\sum_{j=0}^{n-1}\chi_{B}\circ
T^{-j} - \mu(B)\big{|} d\mu = 0.
\]
($\chi_{B}$, the characteristic function of the set  $B$.).  A transformation $T$ is 
\textbf{totally ergodic} when for any $\ell\in\mathbb{N}^{+}$, the 
transformation $T^{\ell}$ is ergodic.  (We use the notation $\mathbb{N} = \{0,1,2,\ldots\}$ for
the natural numbers; $\mathbb{N}^{+} = \{1,2,\ldots\}$ for the positive natural numbers;
and $\mathbb{Z}_{N} = \{0,1,\ldots,N-1\}$ for the $N$-element additive group).
\end{paragr}

\begin{paragr}\textbf{Ergodic Sequences.}
The term \textbf{sequence} shall mean sequence 
in $\mathbb{N}$ that is strictly increasing.  A sequence $\{a_{n}\}$ is 
an \textbf{ergodic sequence} (with respect to a transformation $T$) when 
for all $B\in\mathcal{B}$,
\[
\lim_{n\to\infty} \int\big{|}\frac{1}{n}\sum_{j=0}^{n-1}\chi_{B}\circ
T^{-a_{j}} - \mu(B)\big{|} d\mu = 0.
\]
\end{paragr}

\section{Dynamical Sequences}
\begin{paragr}\textbf{Basic Notions.}
Dynamical sequences were introduced in \cite{CS04}.
\begin{definition}\label{D:dynseq}
A \textbf{\emph{dynamical sequence}}
$\{s_{n,j}\}_{\{r_{n}\}}$ is a doubly-indexed collection of
integers $s_{n,j}$ for $n\in\mathbb{N}$ and $j\in\mathbb{Z}_{r_{n}}$ 
where 
$\{r_{n}\}$ is a given sequence, called the \textbf{\emph{index 
sequence}}, which must have the property that 
$\lim_{n\to\infty}r_{n}=\infty$ (see \cite{CS04}).  The
integer $s_{n,j}$ is the $j^{th}$ \textbf{\emph{element}}
of the dynamical sequence at the $n^{th}$ stage.
\end{definition}
\end{paragr}

\begin{paragr}\textbf{Partial Sums of Dynamical Sequences.}
\begin{notation}
Let $\{s_{n,j}\}_{\{r_{n}\}}$ be a dynamical sequence and 
$n\in\mathbb{N}$, $j\in\mathbb{Z}_{r_{n}}$, $k\in\mathbb{Z}_{r_{n}-j}$.  
The $k^{th}$ partial sum of the $j^{th}$ element at the $n^{th}$ stage is
\[
s_{n,j}^{(k)} = \sum_{z=0}^{k-1}s_{n,j+z}.
\]
\end{notation}
\begin{definition}
Let $\{s_{n,j}\}_{\{r_{n}\}}$ be a dynamical sequence and 
$k\in\mathbb{N}$.  The $k^{th}$
\textbf{\emph{partial sum dynamical sequence}} is the dynamical
sequence $\{s_{n,j}^{(k)}\}_{\{r_{n}-k\}}$
($n$ ``begins'' at the smallest value such that $r_{n}\geq k$) whose
elements are the $k^{th}$ partial sums of
$\{s_{n,j}\}_{\{r_{n}\}}$.
Let $\{k_{n}\}$ be a sequence such that $k_{n} <
r_{n}$ for all $n$.  The $\{k_{n}\}^{th}$
\textbf{\emph{partial sum dynamical sequence}} of
$\{s_{n,j}\}_{\{r_{n}\}}$ is the dynamical sequence $\{s_{n,j}^{(k_{n})}\}_{\{r_{n}-k_{n}\}}$.
\end{definition}
\end{paragr}

\begin{paragr}\textbf{Monotonic Dynamical Sequences.}
\begin{definition}
A dynamical sequence $\{s_{n,j}\}_{\{r_{n}\}}$ is
\textbf{\emph{monotone}} when for every fixed 
$M\in\mathbb{N}$, 
$\lim_{n\to\infty}\frac{1}{r_{n}}\#\{j\in\mathbb{Z}_{r_{n}} : 
\big{|}s_{n,j}\big{|} < M\} = 0$ (the symbol $\#$ denotes cardinality).
\end{definition}
\end{paragr}

\begin{paragr}\textbf{Slicings of Dynamical Sequences.}
\begin{definition}
Let $\{s_{n,j}\}_{\{r_{n}\}}$ be a dynamical sequence and let 
$\{Q_{n}\}$ be a sequence such that 
$\lim_{n\to\infty}\frac{Q_{n}}{r_{n}} = 0$.  
A collection of sets $\Gamma_{n,q} 
\subseteq \mathbb{Z}_{r_{n}}$ and maps $\Psi_{n,q} : \mathbb{Z}_{r_{n}} 
\to \mathbb{Z}$ for $n\in\mathbb{N}$ and $q\in\mathbb{Z}_{Q_{n}}$ is a 
\textbf{\emph{slicing}} of $\{s_{n,j}\}_{\{r_{n}\}}$ when for each 
$n\in\mathbb{N}$ the $\Gamma_{n,q}$ partition $\mathbb{Z}_{r_{n}}$ 
(recall that $\{\Gamma_{q}\}$ partition $\mathbb{Z}_{r}$ when 
$\bigcup\Gamma_{q} = \mathbb{Z}_{r}$ and $\Gamma_{q} \cap 
\Gamma_{q^{\prime}} = \emptyset$ for $q\ne q^{\prime}$), and for each 
$q\in\mathbb{Z}_{Q_{n}}$ there exists $a_{n,q},b_{n,q}\in\mathbb{N}$ 
such 
that for all $j\in\mathbb{Z}_{r_{n}}$, if $a_{n,q} \leq s_{n,j} < 
b_{n,q}$ then $j\in\Gamma_{n,q}$.
\end{definition}
\begin{definition}
Let $\{s_{n,j}\}_{\{r_{n}\}}$ be a dynamical sequence and $\{k_{n}\}$ a 
sequence such that $k_{n} < r_{n}$ for all $n\in\mathbb{N}$.
Let $\{\Gamma_{n,q}\}$ and $\{\Psi_{n,q}\}$ indexed over $\{Q_{n}\}$ be 
a slicing of the partial sum dynamical sequence 
$\{s_{n,j}^{(k_{n})}\}_{\{r_{n} - k_{n}\}}$ and let 
$\{\alpha_{n,q}\}_{\{Q_{n}\}}$ be a 
dynamical sequence such that $\alpha_{n,q} < k_{n}$ for all 
$n\in\mathbb{N}$ and $q\in\mathbb{Z}_{r_{n}-k_{n}}$.  
The dynamical sequences 
$\{s_{n,\Psi_{n,q}(j)}^{(k_{n}-\alpha_{n,q})}\}_{\{\#\Gamma_{n,q}\}}$ 
are an \textbf{\emph{approximate slicing}} of the $\{k_{n}\}^{th}$ partial sum dynamical sequence of $\{s_{n,j}\}_{\{r_{n}\}}$.
\end{definition}
\end{paragr}

\section{Ergodicity on Dynamical Sequences}\label{S:dynavg}

\begin{paragr}\textbf{Ergodic Dynamical Sequences.}
\begin{definition}
A dynamical sequence $\{s_{n,j}\}_{\{r_{n}\}}$ is an 
\textbf{\emph{ergodic dynamical sequence}} (with 
respect to a transformation $T$) when for all $B\in\mathcal{B}$,
\[
\lim_{n\to\infty}\int\big{|}\frac{1}{r_{n}}\sum_{j=0}^{r_{n}-1}\chi_{B}\circ 
T^{-s_{n,j}} - \mu(B)\big{|} d\mu = 0.
\]
\end{definition}
\end{paragr}

\begin{paragr}\textbf{Slice-Ergodicity.}
\begin{definition}
Let $\{s_{n,j}\}_{\{r_{n}\}}$ be a dynamical sequence and $\{k_{n}\}$ be 
a sequence such that $k_{n} < r_{n}$ for all $n\in\mathbb{N}$.  Then
$\{s_{n,j}\}_{\{r_{n}\}}$ is \textbf{\emph{slice-ergodic 
around $\{k_{n}\}$}} (with respect to a 
transformation $T$) if for every approximate slicing of the dynamical 
sequence $\{s_{n,j}^{(k_{n})}\}_{\{r_{n}-k_{n}\}}$ defined by sets  
$\{\Gamma_{n,q}\}$, maps $\{\Psi_{n,q}\}$, sequence $\{Q_{n}\}$, and 
dynamical 
sequence $\{\alpha_{n,q}\}_{\{Q_{n}\}}$, and all $B\in\mathcal{B}$,
\[
\lim_{n\to\infty}\int\big{|}\frac{1}{r_{p_{n}}}\sum_{q=0}^{Q_{n}-1}\sum_{j\in\Gamma_{n,q}}\chi_{B}\circ 
T^{-s_{n,j}^{(k_{n}-\alpha_{n,q})}} - \mu(B)\big{|} d\mu = 0,
\]
i.e., the ergodic average over the approximate slicing tends to zero.
A dynamical sequence is \textbf{\emph{slice-ergodic}} 
when it is slice-ergodic around every 
sequence $\{k_{n}\}$ such that $k_{n} < r_{n}$ for all $n\in\mathbb{N}$.
\end{definition}
\end{paragr}

\begin{paragr}\textbf{Mixing and Ergodic Dynamical Sequences.}
The following  standard generalization of the Blum-Hanson theorem 
from  sequences to dynamical sequences is shown in \cite{CS04}.
\begin{theorem}\label{T:genBH}
A transformation $T$ is mixing if and only if every monotone dynamical 
sequence is ergodic with respect to $T$.
\end{theorem}
\end{paragr}

\section{Rank-One Transformations}\label{S:rankone}

\begin{paragr}\textbf{Cutting and Stacking.}
Begin with $[0,1)$, the only ``level'' in the initial ``column''.  ``Cut'' it into
$r_{0}$ ``sublevels'', pieces of equal length: $[0,\frac{1}{r_{0}})$, 
$[\frac{1}{r_{0}},\frac{2}{r_{0}})$, $\ldots$, 
$[\frac{r_{0}-1}{r_{0}},1)$.
Place $s_{0,0}$ intervals of the same length ``above'' 
$[0,\frac{1}{r_{0}})$, i.e., if $s_{0,0} = 1$ place 
$[1,\frac{r_{0}+1}{r_{0}})$ above
$[0,\frac{1}{r_{0}})$.  Likewise, place $s_{0,j}$ ``spacer'' sublevels 
above each piece.  Now, ``stack'' the resulting subcolumns
from left to right by placing $[0,\frac{1}{r_{0}})$ at the bottom, the 
$s_{0,0}$ spacers above it, $[\frac{1}{r_{0}},\frac{2}{r_{0}})$ above
the topmost of the $s_{0,0}$ spacers and so on, ending with the topmost 
of the $s_{0,r_{0}-1}$ spacers.  This stack of $h_{1} = r_{0} + 
\sum_{j=0}^{r_{0}-1}s_{0,j}$
levels (of length $\frac{1}{r_{0}}$), the second column, defines a map 
$T_{0} : 
[0,1+\frac{1}{r_{0}}\sum_{j=0}^{r_{0}-1}s_{0,j}-\frac{1}{r_{0}}) \to 
[\frac{1}{r_{0}},1+\frac{1}{r_{0}}\sum_{j=0}^{r_{0}-1}s_{0,j})$ that 
sends points directly up one level.

Repeat the process: cut the entire new column into $r_{1}$ subcolumns 
of equal width $\frac{1}{r_{0}r_{1}}$, preserving the stack map on each 
subcolumn;
place $s_{1,j}$ spacers (intervals not yet in the space the same width as the subcolumns) above each subcolumn ($j\in\mathbb{Z}_{r_{1}}$); and stack
the resulting subcolumns from left to right.  Our new column defines a map $T_{1}$ that agrees with $T_{0}$ where it is defined and extends it to all
but the topmost spacer of the rightmost subcolumn.  Iterating this process leads to a transformation $T$ defined on all but a Lebesgue measure zero set.
\end{paragr}

\begin{paragr}\textbf{Construction of Rank-One Transformations.}
A transformation created by \emph{cutting and stacking} as just 
described (with a signle column resulting from each iteration) is a {\bf 
rank-one transformation}.  The reader is referred to 
\cite{Fe97} and \cite{Fr70} for more 
details.  Rank-one transformations are measurable and 
measure-preserving under Lebesgue measure, and are completely defined by a 
dynamical sequence $\{s_{n,j}\}_{\{r_{n}\}}$ where at the $n^{th}$ step we cut
into $r_{n}$ pieces and place $s_{n,j}$ spacers above each subcolumn.  
This $\{s_{n,j}\}_{\{r_{n}\}}$ is the {\bf spacer sequence} for the 
transformation
and $\{r_{n}\}$ is the {\bf cut sequence}.  The {\bf height sequence} 
$\{h_{n}\}$ is the number of levels in each column: $h_{0} = 1$ and $h_{n+1} = r_{n}h_{n} + 
\sum_{j=0}^{r_{n}-1}s_{n,j}$.

We write $I_{n,i}$ to denote the $i^{th}$ level in the $n^{th}$ stack 
($i\in\mathbb{Z}_{h_{n}}$) where $I_{n,0}$ is the bottom level and $T(I_{n,i}) = 
I_{n,i+1}$ and write $C_{n} = \bigcup_{i=0}^{h_{n}-1}I_{n,i}$ to denote the $n^{th}$ column and $S_{n} = C_{n+1} \setminus C_{n}$ to denote the spacers added.
We write $I_{n,i}^{[j]}$ for the $j^{th}$ sublevel of the $i^{th}$ level of the $n^{th}$ column, i.e., $I_{n,0}^{[0]}$ is the leftmost sublevel of the bottom 
level ($I_{n,0}^{[0]} = I_{n+1,0}$ becomes the bottom level of the next column).  Note that $T$ is defined on a finite measure space if and only if 
$\sum_{n=0}^{\infty}\mu(S_{n}) < \infty$ and in that case $T$ is isomorphic to the transformation defined on $[0,1)$ obtained by cutting and stacking in the 
same fashion as $T$ but beginning with $C_{0} = [0,\frac{1}{K})$ where 
$K$ is the measure of the space the original $T$ is defined on.
\end{paragr}

\begin{paragr}\textbf{Rank-One Uniform Mixing.}
Rank-one uniform mixing
involves sums of mixing values over increasingly fine levels.  
Introduced in \cite{CS04}, details and proofs may be found there.
\begin{definition}
Let $T$ be a rank-one transformation with heights $\{h_{n}\}$ and levels 
$\{I_{n,i}\}$ and $\{a_{n}\}$ a sequence.  Set $p_{n}$ such 
that $h_{p_{n}} \leq a_{n} < h_{p_{n}+1}$.  Then
$\{a_{n}\}$ is \textbf{\emph{rank-one uniform mixing}} (with
respect to $T$) when for all $B\in\mathcal{B}$,
\[
\lim_{n\to\infty}\sum_{i=0}^{h_{p_{n}
}-1}\big{|}\mu(T^{a_{n}}(I_{p_{n} ,i}) \cap B) - \mu(I_{p_{n}
,i})\mu(B)\big{|} = 0;
\]
$T$ is \textbf{\emph{rank-one uniform mixing}}
when $\mathbb{N}$ is a rank-one uniform mixing sequence (with respect to 
$T$).
\end{definition}

\begin{proposition}\label{P:unifmix}\emph{\cite{CS04}}
Let $T$ be a rank-one transformation.  If a sequence $\{t_{n}\}$
is rank-one uniform mixing (with respect to $T$) then $\{t_{n}\}$ is 
mixing (with respect to $T$).  Consequently, if $T$ is rank-one uniform 
mixing then $T$ is mixing.
\end{proposition}
\end{paragr}

\begin{paragr}\textbf{Levels of Rank-One Transformations.}
\begin{lemma}\label{L:levels}
Let $T$ be a rank-one transformation with levels 
$\{I_{n,i}\}$, heights $\{h_{n}\}$, and spacers 
$\{s_{n,j}\}_{\{r_{n}\}}$.  Let $p\in\mathbb{N}$, 
$i\in\mathbb{Z}_{h_{p}}$, $j\in\mathbb{Z}_{r_{p}}$, 
$k\in\mathbb{Z}_{r_{p}-j}$, $t\in\mathbb{Z}_{h_{p}-i}$ and $B$ a 
union of levels in $C_{p}$.  Then the following hold:
{\allowdisplaybreaks
\begin{align*}
(i)\quad &I_{p,i} = \bigcup_{j=0}^{r_{p}-1} I_{p,i}^{[j]}; \\
(ii)\quad  &T^{kh_{p} + s_{p,j}^{(k)}}(I_{p,i}^{[j]}) = I_{p,i}^{[j+k]}; \text{and} \\
(iii)\quad &\mu(T^{t}(I_{p,i}^{[j]}) \cap B) = \frac{1}{r_{p}}\mu(T^{t}(I_{p,i}) \cap B).
\end{align*}}
\end{lemma}
\begin{proof}
(i) and (ii) follow from the construction of rank-one 
transformations.  
For (iii), $T^{t}(I_{p,i}^{[j]}) = 
I_{p,i+t}^{[j]}$ and $B \subseteq C_{p}$ so $I_{p,i+t} \subseteq B$ or 
$I_{p,i+t} \cap B = \emptyset$.
\end{proof}

\begin{lemma}\label{L:fulltrick}
For any $p\in\mathbb{N}$, any $\Lambda \subseteq \mathbb{Z}_{h_{N}}$,
any $\Gamma \subseteq \mathbb{N}$, any $B$ a union of levels in $C_{p}$ 
and any maps
$f : \Gamma \to \mathbb{Z}$ and $g : \Gamma \to \mathbb{Z}_{r_{p}}$,
{\allowdisplaybreaks
\begin{align*}
\sum_{i\in\Lambda}\big{|}\sum_{j\in\Gamma} 
&\mu(T^{f(j)}(I_{p,i}^{[g(j)]}) \cap B) - \mu(I_{p,i}^{[g(j)]})
\mu(B)\big{|} \\ &\leq \int \big{|} 
\frac{1}{r_{p}}\sum_{j\in\Gamma}\chi_{B}\circ T^{f(j)} - \mu(B)\big{|}
d\mu + \big{(}\sup_{j\in\Gamma} 
f(j)\big{)}\frac{1}{h_{p}}\frac{\#\Gamma}{r_{p}}.
\end{align*}}
\end{lemma}
\begin{proof}
Lemma \ref{L:levels} (iii) and the definition of integration.
\end{proof}
\end{paragr}

\section{Mixing Sequence Theorems}
\begin{paragr}\textbf{Mixing Height Sequences.}
\begin{theorem}\label{T:mixheights}
Let $T$ be a rank-one transformation with spacers
$\{s_{n,j}\}_{\{r_{n}\}}$ and heights $\{h_{n}\}$ and let 
$k\in\mathbb{N}$.  If 
$\{s_{n,j}^{(k)}\}_{\{r_{n}-k\}}$ is ergodic (with respect to $T$) then 
$\{kh_{n}\}$ is rank-one uniform mixing (with respect to $T$).
\end{theorem}
\begin{proof}
Let $T$, $\{s_{n,j}\}_{\{r_{n}\}}$, $\{h_{n}\}$ and $k$ be as above.  
Let $B$ be a union of levels in $C_{N}$ for some fixed $N\in\mathbb{N}$ 
(note that levels approximate measurable sets).
For any sets $J_{n}\in\mathbb{Z}_{r_{n}}$, 
apply Lemmas \ref{L:levels} (i) then \ref{L:levels} (ii) and 
finally Lemma
\ref{L:fulltrick}, for any $n\geq N$,
{\allowdisplaybreaks
\begin{align*}
\sum_{i=0}^{h_{n}-1}\big{|}&\mu(T^{kh_{n}}(I_{n,i})\cap B) - \mu(I_{n,i})\mu(B)\big{|} \\
&\leq \sum_{i=0}^{h_{n}-1}\big{|}\sum_{j\in\mathbb{Z}_{r_{n}-k}\setminus J_{n}}\mu(T^{kh_{n}}(I_{n,i}^{[j]})\cap B) - \mu(I_{n,i}^{[j]})\mu(B)\big{|} + \frac{k}{r_{n}} + \frac{\# J_{n}}{r_{n}} \\
&= \sum_{i=0}^{h_{n}-1}\big{|}\sum_{j\in\mathbb{Z}_{r_{n}-k}\setminus J_{n}}\mu(T^{s_{n,j}^{(k)}}(I_{n,i}^{[j+k]})\cap B) - \mu(I_{n,i}^{[j+k]})\mu(B)\big{|} + \frac{k}{r_{n}} + \frac{\# J_{n}}{r_{n}} \\
&\leq \int\big{|}\frac{1}{r_{n}}\sum_{j=0}^{r_{n}}\chi_{B}\circ 
T^{s_{n,j}^{(k)}} - \mu(B)\big{|} d\mu + 2\frac{k}{r_{n}} + 2\frac{\# J_{n}}{r_{n}} + \frac{1}{h_{n}}\sup_{j\in\mathbb{Z}_{r_{n}-k}\setminus J_{n}}s_{n,j}^{(k)}.
\end{align*}}

As $k$ is fixed and $\{s_{n,j}^{(k)}\}_{\{r_{n}-k\}}$ is ergodic with 
respect to $T$, we need only show that there exists sets $J_{n} 
\subseteq \mathbb{Z}_{r_{n}}$ such that $\frac{\# J_{n}}{r_{n}} \to 0$ 
and $\frac{1}{h_{n}}\sup_{j\notin J_{n}}s_{n,j}^{(k)} \to 0$.  Suppose 
not.  Then there exists $\delta > 0$ such that $s_{n,j}^{(k)} \geq 
\delta h_{n}$ for at least $\delta r_{n}$ values of $j$ (for infinitely 
many $n$).  But then at least $\frac{1}{k} \delta r_{n}$ values of $j$ 
are such that $s_{n,j} \geq \frac{\delta}{k}$ so $\mu(S_{n}) \geq 
\frac{\delta^{2}}{k^{2}} r_{n} h_{n} \mu(I_{n+1,0}) = 
\frac{\delta^{2}}{k^{2}} \mu(C_{n})$ contradicting that $T$ is defined 
on a finite measure space.
\end{proof}
\end{paragr}

\begin{paragr}\textbf{Mixing Sequences.}
\begin{theorem}\label{T:mixseq}
Let $T$ be a rank-one transformation with spacers
$\{s_{n,j}\}_{\{r_{n}\}}$ and heights $\{h_{n}\}$ and let 
$\{t_{n}\}$ be a sequence.  For each $n\in\mathbb{N}$, set $p_{n}$ 
(uniquely) such that $h_{p_{n}} \leq t_{n} < h_{p_{n}+1}$ and set 
$k_{n}$ (uniquely) such that $k_{n}h_{p_{n}} \leq t_{n} 
< (k_{n}+1)h_{p_{n}}$.  If $\{s_{n,j}\}_{\{r_{n}\}}$ is slice-ergodic 
around $\{k_{n}+1\}$ then the sequence $\{t_{n}\}$ 
is rank-one uniform mixing (with respect to $T$).
\end{theorem}
\begin{proof}
We begin with a brief outline of the method undertaken.  First, 
we dispose of the case when $\frac{k_{n}}{r_{p_{n}}} \to 1$ as it is 
trivial 
from the preceding theorem.  We begin by slicing the 
spacer sequence into blocks of values with difference less than 
$\epsilon_{n}h_{p_{n}}$ ($\epsilon_{n} \to 0$), forming $Q_{n}$ slices 
at 
each stage.  Next, we determine $\alpha_{n,q}$, the number 
of times each subcolumn in the $q^{th}$ block will be mapped through the 
top of the stack under $t_{n}$.  Then we show which sublevel each 
sublevel is mapped to under $T^{t_{n}}$ and or 
each of the three cases arising, and for each $q$, we show that the 
rank-one uniform 
mixing sum is small.  The proof is completed by showing the combined sum 
of the three cases over all the $q$ tends to zero by the 
slice-ergodicity of the spacers.

Let $T$, $\{s_{n,j}\}_{\{r_{n}\}}$, $\{h_{n}\}$, $\{t_{n}\}$, 
$\{p_{n}\}$, and $\{k_{n}\}$ be as above.  Let $\{C_{n}\}$, $\{S_{n}\}$, 
and $\{I_{n,i}\}$ be the columns, spacers, and levels, for $T$, 
respectively.  Let $B$ be a union of levels in $C_{N}$ for some 
$N\in\mathbb{N}$.  Let $n\in\mathbb{N}$ such that $t_{n} 
\geq h_{N}$.  Set $m_{n} = t_{n} - k_{n}h_{p_{n}}$ so 
$m_{n}\in\mathbb{Z}_{h_{p_{n}}}$.  

For $n$ such that $\frac{k_{n}}{r_{p_{n}}} \to 1$, apply Lemma 
\ref{L:levels} (i) and the triangle inequality,
{\allowdisplaybreaks
\begin{align*}
\sum_{i=0}^{h_{p_{n}}-1}&\big{|}\mu(T^{t_{n}}(I_{p_{n},i}) \cap B) - \mu(I_{p_{n},i})\mu(B)\big{|} \\
&\leq \sum_{i=0}^{h_{p_{n}}-1}\sum_{j=0}^{r_{p_{n}}-1}\big{|}\mu(T^{k_{n}h_{p_{n}}+m_{n}}(I_{p_{n},i}^{[j]}) \cap B) - \mu(I_{p_{n},i}^{[j]})\mu(B)\big{|} \\
&\leq \sum_{i=0}^{h_{p_{n}+1}-1}\big{|}\mu(T^{h_{p_{n}+1}}(I_{p_{n}+1,i}) \cap B) - \mu(I_{p_{n}+1,i})\mu(B)\big{|} \\
&\quad\quad + (h_{p_{n}+1} - k_{n}h_{p_{n}} - m_{n})\mu(I_{p_{n}+1,0})
\end{align*}}as there are at most $(h_{p_{n}+1} - k_{n}h_{p_{n}} - m_{n})$ sublevels 
that do not ``map through'' the top of $C_{p_{n}+1}$.  This quantity 
approaches zero as $n\to\infty$ since $\{h_{n}\}$ is rank-one uniform 
mixing 
with respect to $T$ by Theorem \ref{T:mixheights} and since $k_{n} \approx r_{p_{n}}$,
{\allowdisplaybreaks
\begin{align*}
(h_{p_{n}+1} - k_{n}h_{p_{n}} - m_{n})\mu(I_{p_{n}+1,0}) &\leq (h_{p_{n}+1} - k_{n}h_{p_{n}})\mu(I_{p_{n}+1,0}) \\ &\approx (h_{p_{n}+1} - r_{p_{n}}h_{p_{n}})\mu(I_{p_{n}+1,0}) = \mu(S_{p_{n}})
\end{align*}}and $\mu(S_{p_{n}}) \to 0$ as $n\to\infty$ because the final space has finite total measure.

Now consider when $\frac{k_{n}}{r_{p_{n}}}$ is bounded away from 1.  We 
define the sequence $\{\epsilon_{n}\}$ by choosing $\epsilon_{n}\to 0$ as $n
\to \infty$ such that
$\frac{1}{\epsilon_{n}}\frac{\mu(S_{n})}{\mu(C_{n})} \to 0$ as
$n\to\infty$ (possible as
$\frac{\mu(S_{n})}{\mu(C_{n})}\to 0$ since the final space has finite 
total measure).

Let $\Psi_{n} : \mathbb{Z}_{r_{p_{n}}-k_{n}} \to
\mathbb{Z}_{r_{p_{n}}-k_{n}}$ be a map such that
$s_{p_{n},\Psi_{n}(j)}^{(k_{n})} \leq
s_{p_{n},\Psi_{n}(j+1)}^{(k_{n})}$ for all
$j\in\mathbb{Z}_{r_{p_{n}}-k_{n}-1}$.  Set $\ell_{n,0} = 0$ and $\alpha_{n,0} = 1$ and then
proceed inductively to define $\ell_{n,q+1}$ and $\alpha_{n,q+1}$ given $\ell_{n,q}$ and $\alpha_{n,q}$ as
follows: choose $\ell_{n,q+1}$ to be the smallest positive integer
less than $r_{p_{n}}-k_{n}$ such that
\[
s_{p_{n},\Psi_{n}(\ell_{n,q+1})}^{(k_{n}-\alpha_{n,q}+1)} -
s_{p_{n},\Psi_{n}(\ell_{n,q})}^{(k_{n}-\alpha_{n,q})} \geq \epsilon_{n}
h_{p_{n}}
\]
if such an integer exists and choose $\ell_{n,q+1} =
r_{p_{n}}$ and set $Q_{n} = q+1$ if not.   Choose $\alpha_{n,q+1}$ such that
{\allowdisplaybreaks
\begin{align*}
\text{i)}\quad &(\alpha_{n,q}-1)h_{p_{n}} +
s_{p_{n},\Psi_{n}(\ell_{n,q})+k_{n}-\alpha_{n,q}+1}^{(\alpha_{n,q}-1)}
< s_{p_{n},\Psi_{n}(\ell_{n,q})}^{(k_{n})} - m_{n}\text{; and} \\
\text{ii)}\quad &s_{p_{n},\Psi_{n}(\ell_{n,q})}^{(k_{n})} - m_{n} \leq
\alpha_{n,q}h_{p_{n}} +
s_{p_{n},\Psi_{n}(\ell_{n,q})+k_{n}-\alpha_{n,q}}^{(\alpha_{n,q})}.
\end{align*}}Set $\beta_{n,q} = \alpha_{n,q}h_{p_{n}} +
s_{p_{n},\Psi_{n}(\ell_{n,q})+k_{n}-\alpha_{n,q}}^{(\alpha_{n,q})}
- s_{p_{n},\Psi_{n}(\ell_{n,q})}^{(k_{n})} + m_{n}$ and note that
$0 \leq \beta_{n,q} < h_{p_{n}} + s_{p_{n},\Psi_{n}(\ell_{n,q})+k_{n}-\alpha_{n,q}}$.  Set $\beta_{n,q}^{\prime} = \beta_{n,q} - s_{p_{n},\Psi_{n}(\ell_{n,q})+k_{n}-\alpha_{n,q}}$ and note that $\beta_{n,q}^{\prime} < h_{p_{n}}$.

For all
$q\in\mathbb{Z}_{Q_{n}}$, define the sets
\[
\Gamma_{n,q} = \big{\{}j\in\mathbb{Z}_{r_{p_{n}}}
:s_{p_{n},\Psi_{n}(\ell_{n,q})}^{(k_{n})} \leq
s_{p_{n},j}^{(k_{n})} <
s_{p_{n},\Psi_{n}(\ell_{n,q+1})}^{(k_{n})}\big{\}}
\]
and the maps $\Psi_{n,q} : \mathbb{Z}_{\ell_{n,q+1}-\ell_{n,q}}
\to \Gamma_{n,q}$ by $\Psi_{n,q}(j) = \Psi_{n}(j+\ell_{n,q})$ for
all $j\in\mathbb{Z}_{\ell_{n,q+1}-\ell_{n,q}}$.  The resulting
approximate slicing of
$\{s_{p_{n},j}^{(k_{n})}\}_{\{r_{p_{n}}-k_{n}\}}$ are the dynamical 
sequences 
$\{s_{p_{n},\Psi_{n,q}(j)}^{(k_{n})}\}_{\{\ell_{n,q+1}-\ell_{n,q}\}}$
indexed over $q$ by $\{Q_{n}\}$.  

Consider the following using the triangle inequality and Lemma 
\ref{L:levels} (ii):
{\allowdisplaybreaks
\begin{align*}
\sum_{i=0}^{h_{n}-1}&\big{|}\mu(T^{t_{n}}(I_{p_{n},i}) \cap B) - \mu(I_{p_{n},i})\mu(B)\big{|} \\
&= \sum_{i=0}^{h_{n}-1}\big{|}\mu(T^{k_{n}h_{p_{n}}+m_{n}}(I_{p_{n},i}) \cap B) - \mu(I_{p_{n},i})\mu(B)\big{|} \\
&\leq \sum_{i=0}^{h_{p_{n}}-1}\big{|}\sum_{j=0}^{r_{p_{n}}-k_{n}-1}\mu(T^{-s_{p_{n},j}^{(k_{n})}+m_{n}}(I_{p_{n},i}^{[j+k_{n}]}) \cap B) - \mu(I_{p_{n},i}^{[j+k_{n}]})\mu(B)\big{|} \\
&\quad\quad + \sum_{i=0}^{h_{p_{n}}-1}\sum_{j=r_{p_{n}}-k_{n}}^{r_{p_{n}}-1}\big{|}\mu(T^{k_{n}h_{p_{n}}}(I_{p_{n},i}^{[j]}) \cap B) - \mu(I_{p_{n},i}^{[j]})\mu(B)\big{|}.
\end{align*}}The ergodicity of $\{s_{n,j}\}_{\{r_{n}\}}$ with 
respect to $T$ implies that the second summand above tends to zero 
(Theorem \ref{T:mixheights}) using the same argument as above.

Note that
{\allowdisplaybreaks
\begin{align*}
\sum_{i=0}^{h_{p_{n}}-1}&\big{|}\sum_{j=0}^{r_{p_{n}}-k_{n}-1}\mu(T^{-s_{p_{n},j}^{(k_{n})}+m_{n}}(I_{p_{n},i}^{[j+k_{n}]}) \cap B) - \mu(I_{p_{n},i}^{[j+k_{n}]})\mu(B)\big{|} \\
&= \sum_{i=0}^{h_{p_{n}}-1}\big{|}\sum_{q=0}^{Q_{n}-1} \sum_{j\in\Gamma_{n,q}}\mu(T^{-s_{p_{n},j}^{(k_{n})}+m_{n}}(I_{p_{n},i}^{[j+k_{n}]}) \cap B) - \mu(I_{p_{n},i}^{[j+k_{n}]})\mu(B)\big{|} \\
\end{align*}}Now for any $q\in\mathbb{Z}_{Q_{n}}$, any $i\in\mathbb{Z}_{h_{p_{n}}}$ 
and any $j\in\Gamma_{n,q}$, using Lemma \ref{L:levels} (ii),
{\allowdisplaybreaks
\begin{align*}
&T^{-s_{p_{n},j}^{(k_{n})} + m_{n}}(I_{p_{n},i}^{[j+k_{n}]}) \\
&= T^{-\big{(}s_{p_{n},j}^{(k_{n})} - s_{p_{n},\Psi_{n}(\ell_{n,q})}^{(k_{n})}\big{)} - s_{p_{n},\Psi_{n}(\ell_{n,q})}^{(k_{n})} + m_{n}}(I_{p_{n},i}^{[j+k_{n}]}) \\
&= T^{-\big{(}s_{p_{n},j}^{(k_{n})} - s_{p_{n},\Psi_{n}(\ell_{n,q})}^{(k_{n})}\big{)} - s_{p_{n},\Psi_{n}(\ell_{n,q})}^{(k_{n})} + m_{n} + \alpha_{n,q}h_{p_{n}} + s_{p_{n},j+k_{n}-\alpha_{n,q}}^{(\alpha_{n,q})}}(I_{p_{n},i}^{[j+k_{n}-\alpha_{n,q}]}) \\
&= T^{-\big{(}s_{p_{n},j}^{(k_{n})} - s_{p_{n},\Psi_{n}(\ell_{n,q})}^{(k_{n})}\big{)} + \beta_{n,q}} (I_{p_{n},i}^{[j+k_{n}-\alpha_{n,q}]}).
\end{align*}}
If $i < h_{p_{n}} - \beta_{n,q}$ then
\[
T^{-s_{p_{n},j}^{(k_{n})} + m_{n}}(I_{p_{n},i}^{[j+k_{n}]}) = T^{-\big{(}s_{p_{n},j}^{(k_{n}-\alpha_{n,q})} - s_{p_{n},\Psi_{n}(\ell_{n,q})}^{(k_{n}-\alpha_{n,q})}\big{)}}(I_{p_{n},i+\beta_{n,q}}^{[j+k_{n}-\alpha_{n,q}]}).
\]
If $i \geq h_{p_{n}} - \beta_{n,q}^{\prime}$ then
{\allowdisplaybreaks
\begin{align*}
&T^{-s_{p_{n},j}^{(k_{n})} + m_{n}}(I_{p_{n},i}^{[j+k_{n}]}) \\ &= T^{-\big{(}s_{p_{n},j}^{(k_{n}-\alpha_{n,q})} - s_{p_{n},\Psi_{n}(\ell_{n,q})}^{(k_{n}-\alpha_{n,q})}\big{)}+\beta_{n,q}-h_{p_{n}} - s_{p_{n},j+k_{n}-\alpha_{n,q}}}(I_{p_{n},i}^{[j+k_{n}-\alpha_{n,q}+1]}) \\
&= T^{-\big{(}s_{p_{n},j}^{(k_{n}-\alpha_{n,q}+1)} - s_{p_{n},\Psi_{n}(\ell_{n,q})}^{(k_{n}-\alpha_{n,q}+1)}\big{)}+\beta_{n,q}-h_{p_{n}} - s_{p_{n},\Psi_{n}(\ell_{n,q})+k_{n}-\alpha_{n,q}}}(I_{p_{n},i}^{[j+k_{n}-\alpha_{n,q}+1]}) \\
&= T^{-\big{(}s_{p_{n},j}^{(k_{n}-\alpha_{n,q}+1)} - s_{p_{n},\Psi_{n}(\ell_{n,q})}^{(k_{n}-\alpha_{n,q}+1)}\big{)}}(I_{p_{n},i+\beta_{n,q}^{\prime}-h_{p_{n}}}^{[j+k_{n}-\alpha_{n,q}+1]})
\end{align*}}as $i +\beta_{n,q}-h_{p_{n}}^{\prime} \geq 0$ because $i \geq h_{p_{n}} - \beta_{n,q}^{\prime}$ and $i +\beta_{n,q}^{\prime}-h_{p_{n}} < i < h_{p_{n}}$ because $\beta_{n,q}^{\prime} < h_{p_{n}}$.

If $h_{p_{n}} - \beta_{n,q} \leq i < h_{p_{n}} - \beta_{n,q}^{\prime}$ 
then, as above,
{\allowdisplaybreaks
\begin{align*}
&T^{-s_{p_{n},j}^{(k_{n})} + m_{n}}(I_{p_{n},i}^{[j+k_{n}]}) \\ &= T^{-\big{(}s_{p_{n},j}^{(k_{n}-\alpha_{n,q})} - s_{p_{n},\Psi_{n}(\ell_{n,q})}^{(k_{n}-\alpha_{n,q})}\big{)}+\beta_{n,q}-h_{p_{n}} - s_{p_{n},j+k_{n}-\alpha_{n,q}}}(I_{p_{n},i}^{[j+k_{n}-\alpha_{n,q}+1]}) \\
&= T^{-\big{(}s_{p_{n},j}^{(k_{n}-\alpha_{n,q}+1)} - 
s_{p_{n},\Psi_{n}(\ell_{n,q})}^{(k_{n}-\alpha_{n,q})}\big{)}}(I_{p_{n},i+\beta_{n,q}-h_{p_{n}}}^{[j+k_{n}-\alpha_{n,q}+1]}).
\end{align*}}Applying the first case and Lemma \ref{L:fulltrick}, then 
that $T$ is measure-preserving,
{\allowdisplaybreaks
\begin{align*}
\sum_{i=0}^{h_{p_{n}}-\beta_{n,q}}&\big{|}\sum_{q=0}^{Q_{n}-1}\sum_{j\in\Gamma_{n,q}}\mu(T^{-s_{p_{n},j}^{(k_{n})}+m_{n}}(I_{p_{n},i}^{[j+k_{n}]}) 
\cap B) - \mu(I_{p_{n},i}^{[j+k_{n}]})\mu(B)\big{|} \\
&\leq \int
\big{|}\frac{1}{r_{p_{n}}}\sum_{q=0}^{Q_{n}-1}\sum_{j\in\Gamma_{n,q}}\chi_{B}\circ
T^{-\big{(}s_{p_{n},j}^{(k_{n}-\alpha_{n,q})} -
s_{p_{n},\Psi_{n}(\ell_{n,q})}^{(k_{n}-\alpha_{n,q})}\big{)}} -
\mu(B)\big{|} d\mu \\ &\quad\quad\quad + 
\sum_{q=0}^{Q_{n}-1}\big{(}\sup_{j\in\Gamma_{n,q}}
s_{p_{n},j}^{(k_{n}-\alpha_{n,q})} -
s_{p_{n},\Psi_{n}(\ell_{n,q})}^{(k_{n}-\alpha_{n,q})}\big{)}
\frac{1}{h_{p_{n}}}\frac{\#\Gamma_{n,q}}{r_{p_{n}}} \\
&\leq \int 
\big{|}\frac{1}{r_{p_{n}}}\sum_{q=0}^{Q_{n}-1}\sum_{j\in\Gamma_{n,q}}\chi_{B}\circ 
T^{-s_{p_{n},j}^{(k_{n}-\alpha_{n,q})}} - \mu(B)\big{|} d\mu + \epsilon_{p_{n}}.
\end{align*}}
Similarly, for the second and third cases above, we have that
{\allowdisplaybreaks
\begin{align*}
\sum_{i=h_{p_{n}}-\beta_{n,q}^{\prime}}^{h_{p_{n}}-1}&\big{|}\sum_{q=0}^{Q_{n}-1}\sum_{j\in\Gamma_{n,q}}\mu(T^{-s_{p_{n},j}^{(k_{n})}+m_{n}}(I_{p_{n},i}^{[j+k_{n}]}) 
\cap B) - \mu(I_{p_{n},i}^{[j+k_{n}]})\mu(B)\big{|} \\
&\leq \int 
\big{|}\frac{1}{r_{p_{n}}}\sum_{q=0}^{Q_{n}-1}\sum_{j\in\Gamma_{n,q}}\chi_{B}\circ 
T^{-s_{p_{n},j}^{(k_{n}-\alpha_{n,q}+1)}} - \mu(B)\big{|} d\mu + \epsilon_{p_{n}}
\end{align*}}
and
{\allowdisplaybreaks
\begin{align*}
\sum_{i=h_{p_{n}}-\beta_{n,q}}^{h_{p_{n}}-\beta_{n,q}^{\prime} - 
1}&\big{|}\sum_{q=0}^{Q_{n}-1}\sum_{j\in\Gamma_{n,q}}\mu(T^{-s_{p_{n},j}^{(k_{n})}+m_{n}}(I_{p_{n},i}^{[j+k_{n}]}) 
\cap B) - \mu(I_{p_{n},i}^{[j+k_{n}]})\mu(B)\big{|} \\
&\leq 
\int\big{|}\frac{1}{r_{p_{n}}}\sum_{q=0}^{Q_{n}-1}\sum_{j\in\Gamma_{n,q}}\chi_{B}\circ 
T^{-s_{p_{n},j}^{(k_{n}-\alpha_{n,q}+1)}} - \mu(B)\big{|} d\mu + \epsilon_{p_{n}}
\end{align*}}
Combining these three cases, we have that
{\allowdisplaybreaks
\begin{align*}
\sum_{i=0}^{h_{p_{n}}-1}&\big{|}\sum_{q=0}^{Q_{n}-1}\sum_{j\in\Gamma_{n,q}}\mu(T^{-s_{p_{n},j}^{(k_{n})}+m_{n}}(I_{p_{n},i}^{[j+k_{n}]}) 
\cap B) - \mu(I_{p_{n},i}^{[j+k_{n}]})\mu(B)\big{|} \\
&\leq 
\int\big{|}\frac{1}{r_{p_{n}}}\sum_{q=0}^{Q_{n}-1}\sum_{j\in\Gamma_{n,q}}\chi_{B}\circ 
T^{-s_{p_{n},j}^{(k_{n}-\alpha_{n,q})}} - \mu(B)\big{|} d\mu \\
&\quad\quad +
2 \int 
\big{|}\frac{1}{r_{p_{n}}}\sum_{q=0}^{Q_{n}-1}\sum_{j\in\Gamma_{n,q}}\chi_{B}\circ 
T^{-s_{p_{n},j}^{(k_{n}-\alpha_{n,q}+1)}} - \mu(B)\big{|} d\mu + 3\epsilon_{p_{n}}
\end{align*}}

Note that $Q_{n} \leq 
\frac{s_{p_{n},\Psi_{n}(r_{p_{n}}-k_{n}-1)}^{(k_{n})}}{\epsilon_{p_{n}}h_{p_{n}}} 
\leq \frac{s_{p_{n},0}^{(r_{p_{n}})}}{\epsilon_{p_{n}}h_{p_{n}}} = \frac{r_{p_{n}}\mu(S_{p_{n}})}{\epsilon_{p_{n}}\mu(C_{p_{n}})}$
so $\frac{Q_{n}}{r_{p_{n}}} \to 0$ as $n\to\infty$ by the construction 
of $\epsilon_{n}$.  Then the quantities above approach zero by 
slice-ergodicity around $\{k_{n}+1\}$.
\end{proof}
\end{paragr}

\section{Mixing Theorem}
\begin{paragr}\textbf{Mixing Rank-One Transformations.}
Our main theorem lifts the ``restricted growth'' condition from the main 
theorem (Theorem 6) of \cite{CS04}, generalizing it to all rank-one 
transformations.
\begin{theorem}\label{T:main}
For a rank-one transformation $T$, the following are equivalent:
{\allowdisplaybreaks
\begin{align*}
(i)\quad &\text{ $T$ is a mixing transformation;} \\
(ii)\quad  &\text{ $T$ is a rank-one uniform mixing transformation; and} \\
(iii)\quad &\text{ the spacer sequence for $T$  is slice-ergodic (with 
respect to $T$).}
\end{align*}}
\end{theorem}
\begin{proof}
Let $T$ be as above.  If (iii) holds, then Theorem \ref{T:mixseq}
implies that every sequence is rank-one uniform mixing with respect to $T$ so 
(ii) holds.  If (ii) holds, then Proposition \ref{P:unifmix} implies
that i) holds.  Assume that (i) holds
but suppose that (iii) does not.  Let $\{\Gamma_{n,q}\}$ define an 
approximate slicing not ergodic with respect to $T$.  As 
$\frac{Q_{n}}{r_{n}}\to 0$, there exists a sequence $\{q_{n}\}$ such 
that $\#\Gamma_{n,q_{n}} \to \infty$ and so an approximate slice of the 
spacer sequence is monotone but not ergodic with respect to $T$.  
Theorem \ref{T:genBH} then yields a contradiction.
\end{proof}
\end{paragr}

\section{Power Ergodicity}

\begin{paragr}\textbf{Power Ergodicity.}
Power ergodicity is all powers of an ergodic 
transformation being ``uniformly'' ergodic in the sense that the 
ergodic averages converge uniformly to zero.  Earlier results on 
specific rank-one mixing used precursors to this notion, including the 
uniform Ces\`{a}ro property used in \cite{AF92} (and implicitly in 
\cite{Ad98}) and power uniform ergodicity in \cite{CS04}.

\begin{definition}A transformation $T$ is \textbf{\emph{power ergodic}}
when for all $B\in\mathcal{B}$,
\[
\lim_{n\to\infty} \sup_{k\in\mathbb{N}}\int
\big{|}\frac{1}{n}\sum_{j=0}^{n-1}\chi_{B}\circ T^{-jk} -
\mu(B)\big{|} d\mu = 0.
\]
\end{definition}
\end{paragr}

\begin{paragr}\textbf{Weak Power Ergodicity.}
\begin{definition}
A transformation $T$ is \textbf{\emph{weak power ergodic}} when for 
every sequence $\{k_{n}\}$ such that $\lim_{n\to\infty}\frac{k_{n}}{n} < 
\infty$ and all $B\in\mathcal{B}$,
\[
\lim_{n\to\infty}\int\big{|}\frac{1}{n}\sum_{j=0}^{n-1}\chi_{B}\circ 
T^{-jk_{n}} - \mu(B)\big{|} d\mu = 0.
\]
\end{definition}

Weak power ergodicity was introduced in \cite{CS04} as 
``power uniform ergodicity'' however in light of the power ergodic 
property this name is more accurate.

\begin{theorem}\label{T:wkpe}
Let $T$ be a rank-one transformation such that for each fixed 
$k\in\mathbb{N}$ the $k^{th}$ partial sum sequence of the spacer 
sequence is ergodic with respect to $T$.  Then $T$ is weak power 
ergodic.
\end{theorem}
\begin{proof}
Let $T$ be a rank-one transformation with spacer sequence 
$\{s_{n,j}\}_{\{r_{n}\}}$ and height sequence $\{h_{n}\}$ and let 
$k\in\mathbb{N}$.  Note that 
$\frac{1}{h_{n}}s_{n,j}^{(k)} \to 0$ as $n\to\infty$ for a density one 
set of $j \in \mathbb{Z}_{r_{n}-k}$ because $T$ is defined on a finite 
measure space (details are left to the reader).  Applying Proposition 
7.3 of \cite{CS04} to those $j$ (and ignoring the rest, a zero measure 
set) yields the conclusion.
\end{proof}
\end{paragr}

\begin{paragr}\textbf{Power Ergodicity Theorem.}
\begin{theorem}\label{T:pe}
Let $T$ be a rank-one transformation with spacers
$\{s_{n,j}\}_{\{r_{n}\}}$ such that for any sequence $\{k_{n}\}$ where 
$\lim_{n\to\infty}\frac{k_{n}}{r_{n}} = 0$, the partial sum 
dynamical sequence $\{s_{n,j}^{(k_{n})}\}_{\{r_{n}-k_{n}\}}$ is ergodic 
with 
respect to $T$.  Then $T$ is power ergodic.
\end{theorem}

\begin{lemma}\label{L:pLR}\emph{\cite{Ad98}} \textbf{\emph{(Block Lemma)}}
Let $T$ be a measure-preserving transformation and $B\in\mathcal{B}$.  
Then for any $R, L, p\in\mathbb{N}$,
\[
\int \big{|}\frac{1}{R}\sum_{j=0}^{R-1}\chi_{B}\circ T^{-j} -
\mu(B)\big{|} d\mu \leq \int
\big{|}\frac{1}{L}\sum_{j=0}^{L-1}\chi_{B}\circ T^{-jp} -
\mu(B)\big{|} d\mu + \frac{pL}{R}.
\]
\end{lemma}

\begin{proof} (of Theorem \ref{T:pe})
Let $T$ be a rank-one transformation with spacer sequence 
$\{s_{n,j}\}_{\{r_{n}\}}$, let $\{k_{n}\}$ be an arbitrary sequence and 
$B\in\mathcal{B}$.  
For each $n\in\mathbb{N}$, set $p_{n}, q_{n}, x_{n} \in\mathbb{N}$ such 
that $h_{p_{n}} < k_{n} \leq h_{p_{n}+1}$, $h_{p_{n}+1} \leq k_{n}q_{n} 
< 2 h_{p_{n}+1}$ and $x_{n}h_{p_{n}} \leq k_{n} < (x_{n}+1)h_{p_{n}}$.  

First, consider the case when $\frac{q_{n}}{n} \to 0$.  Fix $\epsilon > 
0$ and choose $L, N_{0} \in \mathbb{N}$ such that $\frac{1}{L} < 
\epsilon$ and $\frac{q_{n}L}{n} < \epsilon$ for $n \geq N_{0}$.  For 
each fixed $\ell < 2L$, the sequence $\{\ell h_{p_{n}+1}\}$ is mixing by 
Theorem \ref{T:mixheights} as $\{s_{n,j}^{(\ell)}\}_{\{r_{n}-\ell\}}$ is 
ergodic with respect to $T$ by hypothesis.  Since $\ell h_{p_{n}+1} \leq 
\ell k_{n}q_{n} < 2\ell h_{p_{n}+1}$, the sequence $\{\ell k_{n}q_{n}\}$ 
is then mixing (following from the construction of rank-one 
transformations).  Let $N\in\mathbb{N}$ such that 
$\big{|}\mu(T^{\ell 
k_{n}q_{n}}(B) \cap B) - \mu(B)\mu(B)\big{|} < \epsilon$ for all $0 <
\ell < L$.  Then, for $n \geq \max(N_{0},N)$, applying the 
Block Lemma (Lemma \ref{L:pLR}) and the H\"{o}lder Inequality,
{\allowdisplaybreaks
\begin{align*}
\int&\big{|}\frac{1}{n}\sum_{j=0}^{n-1}\chi_{B}\circ T^{-jk_{n}} - 
\mu(B)\big{|} d\mu \\ &\leq 
\int\big{|}\frac{1}{L}\sum_{\ell=0}^{L-1}\chi_{B}\circ T^{-\ell 
k_{n}q_{n}} - \mu(B)\big{|} d\mu + \frac{q_{n}L}{n} \\
&\leq 
\Big{[}\frac{1}{L}\sum_{\ell=-L+1}^{L-1}\frac{L-\ell}{L}\big{(}\mu(T^{\ell 
k_{n}q_{n}}(B) \cap B) - \mu(B)\mu(B)\big{)}\Big{]}^{\frac{1}{2}} + 
\epsilon < 
2\epsilon.
\end{align*}}

Now, consider the case when $\frac{q_{n}}{n} \geq \delta$ for some 
$\delta > 0$.  Since $\frac{q_{n}}{n} < \frac{2h_{p_{n}+1}}{nk_{n}} \leq 
\frac{2h_{p_{n}+1}}{nx_{n}h_{p_{n}}} \approx \frac{2r_{p_{n}}}{nx_{n}}$ 
(by finite measure-preserving), $\frac{x_{n}}{r_{p_{n}}} < 
\frac{1}{\delta n}$.  Choose a sequence $\{L_{n}\}$ such that 
$L_{n}\to\infty$ and $\frac{L_{n}}{n}\to 0$.  For any sequence 
$\{\ell_{n}\}$ such that $\ell_{n} < L_{n}$, we see that 
$\frac{\ell_{n}x_{n}}{r_{p_{n}}} < \frac{L_{n}}{n\delta} \to 0$.  By 
hypothesis, this means that 
$\{s_{p_{n},j}^{(\ell_{n}x_{n})}\}_{\{r_{p_{n}} - \ell_{n}x_{n}\}}$ is 
ergodic with respect to $T$.  Theorem 5 of \cite{CS04} then yields that 
$\{\ell_{n}k_{n}\}$ is a mixing sequence.  Applying the Block Lemma 
(Lemma \ref{L:pLR}) and the H\"{o}lder Inequality,
{\allowdisplaybreaks
\begin{align*}
\int&\big{|}\frac{1}{n}\sum_{j=0}^{n-1}\chi_{B}\circ T^{-jk_{n}} - 
\mu(B)\big{|} d\mu \\ &\leq 
\int\big{|}\frac{1}{L_{n}}\sum_{\ell=0}^{L_{n}-1}\chi_{B}\circ T^{-\ell 
k_{n}} - \mu(B)\big{|} d\mu + \frac{L_{n}}{n} \\
&\leq 
\Big{[}\frac{1}{L_{n}}\sum_{\ell=-L_{n}+1}^{L_{n}-1}\frac{L_{n}-\ell}{L_{n}}\big{(}\mu(T^{\ell 
k_{n}}(B) \cap B) - \mu(B)\mu(B)\big{)}\Big{]}^{\frac{1}{2}} + 
\frac{L_{n}}{n} \to 0.
\end{align*}}
\end{proof}
\end{paragr}

\begin{paragr}\textbf{Polynomial Power Ergodicity}
Polynomial power ergodicity is the ``polynomial 
powers'' of a transformation being ``uniformly ergodic''.
The term \textbf{polynomial} shall mean polynomials with 
rational coefficients that map integers to integers.
\begin{definition}A transformation $T$ is \textbf{\emph{polynomial power 
ergodic}} when for all
sequences of polynomials $\{p_{n}\}$ of bounded degree and 
all $B\in\mathcal{B}$,
\[
\lim_{n\to\infty} \int
\big{|}\frac{1}{n}\sum_{j=0}^{n-1}\chi_{B}\circ T^{-p_{n}(j)} -
\mu(B)\big{|} d\mu = 0;
\]
$T$ is \textbf{\emph{weak polynomial power ergodic}} when the polynomial 
power ergodicity condition holds for polynomial sequences $\{p_{n}\}$ 
of bounded degree such that $\lim_{n\to\infty}\frac{c_{n}}{n} < \infty$ 
where $\{c_{n}\}$ are the lead coefficients of the $\{p_{n}\}$.
\end{definition}

\begin{theorem}\emph{\cite{CS04}}\label{T:wkpolype}
Let $T$ be a transformation that is weak power ergodic.  Then
$T$ is weak polynomial power ergodic.
\end{theorem}
\begin{theorem}\label{T:polype}
Let $T$ be a transformation that is power ergodic.  Then $T$ is 
polynomial power ergodic.
\end{theorem}
\begin{proof}Identical to that of Theorem 8 in \cite{CS04} (which is 
Theorem \ref{T:wkpolype} above)  omitting the condition for weakness.
\end{proof}
\end{paragr}

\section{Mixing Rank-One Transformations}
\label{S:MixRankOne}
\begin{paragr}\textbf{Staircase Transformations.}
Let $\{r_{n}\}$ be a sequence and $T$ a rank-one transformation with 
spacer sequence $\{s_{n,j}\}_{\{r_{n}\}}$ given by $s_{n,j} = j$.  Then 
$T$ is a \textbf{staircase transformation}.
\end{paragr}
\begin{paragr}\textbf{Polynomial Staircase Transformations.}
Let $\{p_{n}\}$ be a sequence of polynomials with
bounded degree.  A rank-one transformation with spacer sequence 
$\{s_{n,j}\}_{\{r_{n}\}}$ given by $s_{n,j} = p_{n}(j)$ is a 
\textbf{polynomial staircase transformation}.  We require that the 
polynomials be such that for every $L\in\mathbb{N}$, 
\[\limsup_{n\to\infty}\frac{1}{r_{n}}\#\{j\in\mathbb{Z}_{r_{n}} : L 
\quad\text{divides}\quad p_{n}(j+1) - p_{n}(j)\} < 1\] and such that 
$\lim_{n\to\infty}\frac{c_{n}}{n} < \infty$, where $\{c_{n}\}$ are the 
lead coefficients of the polynomials.

\begin{remark}
Adams and Friedman in \cite{AF92} 
introduced   polynomial staircase transformations as models for rank-one mixing transformations.   A
class of mixing polynomial staircase transformations 
was  studied in the first author's undergraduate thesis at Williams College
in 2001.  Polynomial staircase transformations with restricted growth
were shown to be mixing in \cite{CS04}.   After  \cite{CS04} was  accepted,  
 Ryzhikov wrote the authors 
that a similar class 
of regular behavior  polynomial staircases 
 was shown to be mixing by him 
in 2002,  but a preprint was not available \cite{Ry03}.

\end{remark}

%\begin{remark}
%Staircase transformations are linear polynomial staircase transformations.
%\end{remark}
The following theorem is essentially the polynomial ergodic theorem of 
Furstenberg \cite[p. 70]{Fu81} applied to a sequence of polynomials with 
uniformly bounded coefficients (see \cite{CS04} for a short proof).
\begin{theorem}\label{T:furst}\emph{[Furstenberg]}
Let $\{p_{n}\}$ be a sequence of polynomials of 
bounded degree and uniformly bounded coefficients.  Then a dynamical 
sequence $\{s_{n,j}\}_{\{r_{n}\}}$ given by $s_{n,j} = p_{n}(j)$ is 
ergodic with 
respect to any totally ergodic transformation.
\end{theorem}
\begin{theorem}\label{T:polyR1mix}
Polynomial staircase transformations are mixing.
\end{theorem}
\begin{proof}
Let $T$ be a polynomial staircase transformation.  Then $T$ is totally 
ergodic by Theorem 7 of \cite{CS04}.  Let $\{p_{n}\}$ be the polynomials 
defining the spacer sequence $\{s_{n,j}\}_{\{r_{n}\}}$ of degree at most 
$D\in\mathbb{N}$ and let $\{c_{n,a}\}$ for $a\in\mathbb{Z}_{D+1}$ be 
the coefficients of the $\{p_{n}\}$.  Then, for any 
$j\in\mathbb{Z}_{r_{n}}, 
k\in\mathbb{Z}_{r_{n}-j}$,
{\allowdisplaybreaks
\begin{align*}
s_{n,j}^{(k)} &= \sum_{z=0}^{k-1}p_{n}(j+z) = 
\sum_{z=0}^{k-1}\sum_{a=0}^{D}c_{n,a}(j+z)^{a} =
\sum_{z=0}^{k-1}\sum_{a=0}^{D}\sum_{b=0}^{a} 
c_{n,a}\binom{a}{b}j^{b}z^{a-b} \\ &= 
\sum_{b=0}^{D}\big{(}\sum_{z=0}^{k-1}\sum_{a=b}^{D}c_{n,a}\binom{a}{b}z^{a-b}\big{)}j^{b} 
= p_{n,k}(j)
\end{align*}}are polynomials of degree at most $D$ in $j$ with 
lead coefficients $k c_{n,D}$.

$T$ is weak power ergodic by Theorem \ref{T:wkpe} as for each fixed 
$k\in\mathbb{N}$ the $p_{n,k}$ are ergodic (Theorem \ref{T:furst}) 
since the coefficients of $p_{n,k}$ are uniformly bounded when $k$ is 
fixed.  
$T$ is then weak polynomial power ergodic by Theorem \ref{T:wkpolype}.  
For any sequence $\{k_{n}\}$ such that $\frac{k_{n}}{r_{n}} \to 0$, the 
partial sum sequence $\{s_{n,j}^{(k_{n})}\}_{\{r_{n} - k_{n}\}}$ is 
ergodic with respect to $T$ by weak polynomial power ergodicity (as 
$s_{n,j}^{(k_{n})} = p_{n,k_{n}}(j)$ where $p_{n,k_{n}}$ have lead 
coefficients $k_{n}c_{n}$ and $\frac{k_{n}c_{n}}{n} \to 0$ as 
$\lim_{n\to\infty}\frac{c_{n}}{n} < \infty$).  Hence, $T$ is power 
ergodic by Theorem \ref{T:pe} and so is polynomial power ergodic by 
Theorem \ref{T:polype}.

For any $\{k_{n}\}$ and approximate slicing 
$\{\Gamma_{n,q}\}, \{\Psi_{n,q}\}, 
\{\alpha_{n,q}\}_{\{Q_{n}\}}$ of the spacer sequence around $\{k_{n}\}$, 
each slice 
$\{s_{n,\Psi_{n,q}(j)}^{(k_{n}-\alpha_{n,q})}\}_{\{\#\Gamma_{n,q}\}}$ is 
itself a polynomial sequence (details are left to the reader) of degree 
at most $D$ .  Since 
$\frac{Q_{n}}{r_{n}}\to 0$, $\#\Gamma_{n,q} \to \infty$ uniformly over a 
density one set of $q$, the slices then uniformly tend to zero by 
polynomial power ergodicity.  Theorem \ref{T:main} then yields the 
result.
\end{proof}
\end{paragr}

 As mentioned in the introduction, Ryzhikov
wrote the authors  that  he had proved the following Corollary~\ref{C:adamstairs} 
 \cite{Ry03}; it answers the question asked by
Adams in \cite{Ad98}.

\begin{corollary}\label{C:adamstairs}
Staircase transformations are mixing.
\end{corollary}

\section{Specific Mixing Transformations}
\begin{paragr}\textbf{Criterion for Finite Measure on Rank-One Transformations}
\begin{proposition}\label{P:finmeas}
A rank-one transformation with spacer sequence 
$\{s_{n,j}\}_{\{r_{n}\}}$ and heights $\{h_{n}\}$ is defined on
a finite measure space if and only if
\[
\sum_{n=0}^{\infty}\frac{\bar{s}_{n}}{h_{n}} < \infty 
\quad\big{(}\text{where}\quad \bar{s}_{n} = 
\frac{1}{r_{n}}\sum_{j=0}^{r_{n}-1}s_{n,j}\big{)}.
\]
\end{proposition}
\begin{proof}
Let $T$, $\{s_{n,j}\}_{\{r_{n}\}}$, $\{\bar{s}_{n}\}$, and $\{h_{n}\}$ 
be as above and let $(X,\mu)$ be the space $T$ is defined on.  Let 
$\{C_{n}\}$
denote the columns of the construction as sets, $\{I_{n}\}$ the base levels of the columns, and
$\{S_{n}\}$ the spacers added (so $S_{n} = C_{n+1} \setminus C_{n}$).  We
see that
\[
\mu(S_{n}) = \sum_{j=0}^{r_{n}-1}s_{n,j}\mu(I_{n+1}) = \big{(}\frac{1}{r_{n}}\sum_{j=0}^{r_{n}-1}s_{n,j}\big{)}\mu(I_{n}) = \frac{\bar{s}_{n}}{h_{n}}\mu(C_{n})
\]
and so
\[
\frac{\mu(C_{n+1})}{\mu(C_{n})} = \frac{\mu(C_{n})+\mu(S_{n})}{\mu(C_{n})} = 1+\frac{\bar{s}_{n}}{h_{n}}.
\]
Then,
\[
\log\Big{(}\frac{\mu(X)}{\mu(C_{0})}\Big{)} =
\log\big{(}\prod_{n=0}^{\infty}\frac{\mu(C_{n+1})}{\mu(C_{n})}\Big{)} =
\sum_{n=0}^{\infty}\log(1+\frac{\bar{s}_{n}}{h_{n}}) \approx \sum_{n=0}^{\infty}\frac{\bar{s}_{n}}{h_{n}}
\]
using the approximation $\log(1+\epsilon)\approx\epsilon$ for small $\epsilon$.

Since $\mu(X) < \infty$ if and only if $\frac{\mu(X)}{\mu(C_{0})}<\infty$, the result follows.
\end{proof}
\end{paragr}

\begin{paragr}\textbf{Specific Examples of Mixing Transformations.}
\begin{definition}
Let $D\in\mathbb{N}$ and $\delta\in\mathbb{R}^{+}$.  The rank-one 
transformation $T_{D,\delta}$ with spacers 
$\{s_{n,j}\}_{\{r_{n}\}}$ given by $s_{n,j} = 
j^{D}$ and $r_{n} = 
\big{\lfloor}h_{n}^{\frac{1}{D+\delta}}\big{\rfloor}$
(where $\{h_{n}\}$ is the heights for $T_{D,\delta}$) is a 
\textbf{\emph{simple polynomial staircase transformation}}.
\end{definition}

\begin{theorem}
Simple polynomial staircase transformations are mixing.
\end{theorem}
\begin{proof}
By Theorem \ref{T:polyR1mix}, we need only show the transformations are 
defined on a 
finite measure space.  Let $D\in\mathbb{N}$ and 
$\delta\in\mathbb{R}^{+}$.  Let $\{s_{n,j}\}_{\{r_{n}\}}$ be the 
spacers and $\{h_{n}\}$ the heights for $T_{D,\delta}$.  Now,
$\sum_{j=0}^{r_{n}-1}s_{n,j} = \sum_{j=0}^{r_{n}-1}j^{D} \approx 
r_{n}^{D+1}$.
Since $h_{n+1} = r_{n}h_{n} + \sum_{j=0}^{r_{n}-1}s_{n,j}$, we see that
$h_{n} \geq \prod_{z=0}^{n-1}r_{z} \geq 2^{n}$.  Then,
\[
\sum_{n=0}^{\infty}\frac{\bar{s}_{n}}{h_{n}} \approx 
\sum_{n=0}^{\infty}\frac{1}{r_{n}h_{n}}r_{n}^{D+1} = \sum_{n=0}^{\infty}\frac{r_{n}^{D}}{h_{n}} 
\approx \sum_{n=0}^{\infty}h_{n}^{\frac{D}{D+\delta}-1} \leq 
\sum_{n=0}^{\infty}\big{(}2^{\frac{\delta}{D+\delta}}\big{)}^{-n} < \infty
\]
by the convergence of geometric series.  Proposition \ref{P:finmeas} 
completes the proof.
\end{proof}
\end{paragr}

\begin{paragr}\textbf{Ornstein's Transformation.}
Ornstein's original construction of rank-one mixing transformations 
involved placing spacer levels randomly using a uniform distribution so 
that almost surely the resulting transformation is mixing.  The uniform 
distribution can be equally well interpreted as meaning that the spacer 
sequence is almost surely slice-ergodic so that mixing for these 
transformations follows from our theorem.  The reader is referred to 
\cite{CS04} for details.
\end{paragr}

\begin{paragr}\textbf{A Note on Restricted Growth.}
The restricted growth condition in \cite{CS04} is equivalent to 
$\frac{r_{n}\bar{s}_{n}}{h_{n}}\to 0$ for polynomial staircase 
transformations 
(details are left to the reader).  For $T_{D,\delta}$, we see that
$\frac{r_{n}\bar{s}_{n}}{h_{n}} \approx \frac{r_{n}^{D+1}}{h_{n}} 
\approx 
h_{n}^{\frac{D+1}{D+\delta}-1} = h_{n}^{\frac{1-\delta}{D+\delta}}$
and so $T_{D,\delta}$ has restricted growth if and only if $\delta > 1$.  
Hence, the theorems of \cite{CS04} apply only to $T_{D,\delta}$ with 
$\delta > 1$.  The $T_{D,\delta}$ for $0 < \delta \leq 1$ are not 
provably mixing using previous results.
In particular, the staircase transformations $T_{1,\delta}$ for $0 < 
\delta \leq 1$ are mixing (the case $\delta > 1$ was first shown by 
Adams in \cite{Ad98}; another proof was given in \cite{CS04}).
\end{paragr}

\bibliography{MixingRankOneFeb6}

\end{document}